\documentclass{nashboro}
\title{Monotonicity Analysis over Chains and Curves}
\author{Dan Kucerovsky and Daniel Lemire}
\def\shortitle{Monotonicity Analysis}
\def\shortauthor{D.~Kucerovsky and D.~Lemire}
\newcommand{\komment}[1]{}

\usepackage{url}
\usepackage{amsmath}\usepackage{amsfonts}
\usepackage{epsfig}

\newenvironment{proof}{\pf}{\eop}
\begin{document}
\maketitle


\def\backsh{{\tt \char'134}}
\def\co{{\cal O}}
\def\conf{Conference~}
\def\cP{{\cal P}}
\def\ct{{\cal T}}
\def\ctil{\tilde{\cal T}}
\def\lbr{\char'173}
\def\ms{{\medskip}}
\def\rbr{\char'175}
\def\oln{\co(n \log(n))}
\def\sm{{\smallskip}}
\def\tb#1{{\tt \backsh #1}}
\def\RR{{{\mathbb R}}}

\begin{abstract}
Chains are vector-valued signals sampling a curve. They  are important to
motion signal processing
and to many scientific applications including location sensors.
We propose a novel measure of smoothness for chains curves
by generalizing the
scalar-valued
concept of monotonicity. Monotonicity can be
defined by the connectedness of the inverse image
of balls. This definition is coordinate-invariant and can be computed efficiently
over chains. Monotone curves may be discontinuous, but continuous
monotone curves are differentiable a.e. 
Over chains,  a simple sphere-preserving filter is shown to never decrease
the degree of monotonicity. It outperforms moving average filters over a
synthetic data set. Applications include Time Series Segmentation,
chain reconstruction from unordered data points, Optical Character Recognition, and Pattern Matching.
\end{abstract}


\section{Introduction}
Monotonicity is one of the simplest property a signal may have. It
offers a powerful qualitative description  (``it
goes up,'' ``it goes down'').
Given data coming in from either sensors or
from a numerical simulation, monotonicity is independent of the sampling
frequency and is robust with
respect to missing data~\cite{LearningFromIncompleteData}. 
Many geometrical objects such as curves are
 typically defined in a parametrization-independent way which makes
monotonicity
  appealing.

In this paper, we are concerned with discretely sampled curves (which we
call chains)
such as  the trajectory of a particle in some vector
space. This problem has applications in motion capture and tracking~\cite{jehee2001,whalenet}.

\begin{figure}[b]
\center{\includegraphics[width=0.60\columnwidth]{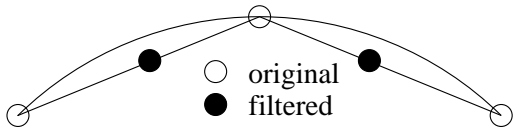}}
\caption{\label{filteredcircle} {\small Given samples on circle, a simple
moving average does not preserve the embedding.}}
\vspace*{-0.5cm}
\end{figure}

We expect a ``smooth'' scalar-valued signal not to change too quickly: it
should be locally constant. Therefore, classical low pass filters such as
the moving average (MA) are often sufficient to help
smooth signals.
Unfortunately, ``smooth'' chains are not locally constant: consider a
loosely sampled circle (see Fig.~\ref{filteredcircle}). 
Moreover, a chain may lie on a sphere or other higher dimensional surface
and we may need to preserve this embedding. In Fig.~\ref{filteredcircle},
a chain on a circle is filtered using
a moving average: we see that the filtered chain can, at best, follow a
circle of a smaller radius. A filter is sphere-preserving (resp.
circle-preserving) if, when the input data points are on a sphere (resp.
circle), the filtered data points also lie on the same sphere (resp.
circle). It is readily shown that no linear filter except the identity can be sphere-preserving
(SP) or circle-preserving (CP). In general, an SP filter is CP.
We offer a simple SP filter in Section~5.



One of the main contribution of this paper is to provide a generalization
of the concept of monotonicity which applies to vector-valued signals
and to curves.
This definition is shown to be robust with respect to removal of data
points and to be efficiently computed. Over curves, we show
that monotone curves have many of the same properties as monotone functions
as far as continuity and differentiability are concerned.
We also propose a SP filter which we show to never decrease the degree of
monotonicity. Experimentally, we show that the degree of monotonicity is
inversely correlated with noise and we compare the SP filter with simple
MA filters, proving the nonlinear SP filter is a good choice when noise
levels are low. Applications of this work include chain reconstruction 
from unordered data points~\cite{289873} and Optical Character Recognition~\cite{201651}.

\section{Related Work}
A motion signal is comprised of
two components: orientation and translation.
The orientation vector indicates
where the object is facing, whereas the translation component determines
the object's location.
  Recent work has focused on smoothing the orientation
vectors~\cite{jehee2001,Hsieh2002}, whereas the results of the present
paper apply equally well to orientation vectors (points on the surface of a unit sphere) 
as to arbitrary translation signals.

In~\cite{monotonevector, pointset,setofpointshakimi},
the authors chose to define monotonicity for curves or chains
with an arbitrary
direction vector: a curve is monotone if its projection on a line is does not backtrack.
While this is a sensible choice given the
lack of definition elsewhere, we argue that not all applications
support an arbitrary direction that can be used to define monotonicity.

The definition of monotonicity has been extended to real-valued
functions~\cite{brooks2004,Carr2000,chiang2003,hong1996,morse69,Kreveld1997}
($f:\mathbb{R}^n \rightarrow
\mathbb{R}$)
by using contour lines (or surfaces) but the idea does not immediately
generalize to curves and chains.

One approach to chain smoothing is to use B-splines and Bezier curves with
the $L^2$ norm~\cite{Smoothingandmatching1994}. Correspondingly, we could
measure the ``smoothness'' of a chain by measuring how closely one can fit
it to a smooth curve. Our approach differs in that we do not use polygonal
approximations or curve fitting: we consider chains to be first-class
citizens.  

\section{Monotone Curves}

Recall that a function $f:\mathbb{R}\rightarrow\mathbb{R}$ is said
to be monotone increasing if $f(x)\geq f(y)$ whenever $x\geq y$
and monotone decreasing if $f(x)\leq f(y)$ whenever $x\geq y$. A
 monotone increasing or monotone decreasing function
is said to be monotone. Recall that
$B=\left\{ x:|x-a|\leq R\right\}$ is called a (closed) ball of radius
$R$ centered around $a$: in the multidimensional case,
the ball is a generalization of the (closed) interval.

\begin{proposition}
$f:\mathbb{R}\rightarrow\mathbb{R}$ is monotone if and only if $f^{-1}(B)$
is connected for all balls $B$.
\end{proposition}
\komment{
\begin{proof}
Suppose $f^{-1}(B)$ is connected for all balls $B$. Then, if $f$ is not
monotone,
there are 3 values $a,b,c$ with $f(a)=f(c)\neq f(b)$ and $a<b<c$.
We can suppose $f$ is uniformly continuous on $[a,b]$. Let
$\epsilon=\frac{1}{2} \vert f(c) = f(a) \vert$, thus there is a $\delta>0$
such
that on $[a,c]$, $\vert f(x)-f(y)\vert < \epsilon$ whenever $\vert x-y
\vert <\delta$,
by uniform continuity. Now, let $B$ be the ball about $f(a)$ of radius
$\epsilon$.
Clearly $f(b)\notin B$ so $b \notin f^{-1}(B)$. But, $f(a)=f(c)$ is in
$B$,
so $a,c\in f^{-1}(B)$. This shows that $f^{-1}(B)$ is a disconnected set.
Conversely, if $f$ is monotone, then clearly $f^{-1}((a,b))$ is an
interval, hence
the inverse image of a ball is connected.
\end{proof}}

An arc-length parametrized curve $s:t\rightarrow s(t)$ is $R$-monotone for $R>0$
if the inverse image of any closed ball of radius at most $R$, under $s$,
is connected. Straight lines are $R$-monotone for all $R>0$. As motivation
the discrete case, we
want to compare monotone curves with monotone functions. Monotone functions $f: \mathbb{R}\rightarrow \mathbb{R}$ are differentiable almost everywhere, and they do not have to be continuous.
 $R$-monotone also do not have to be continuous: the curve
 $s(t)=(f(t),f(t))$ where  $f'(t)=1$ a.e. is $R$-monotone for all $R>0$.
Moreover, they are also differentiable a.e. as the next proposition shows.

\begin{proposition}
Continuous $R$-monotone curves are differentiable a.e.
\end{proposition}
\begin{proof}
Take any point $x$ in the (open) domain of the curve $s$. Choose another point $y$ so that the 
arc-length $y-x$ over $s$ is smaller than $R$. Consider any point $z$ on $s$ between $y$ and $x$,
then $z$ must be contained in all balls of radius $R$ containing both $x$ and $y$. It follows
that $s$ must be differentiable from the left at $x$. Similarly, $s$ is differentiable from the right
at $x$. If the two derivative from the left and from the right do not match, 
then it is possible to find $y$ and $y'$ close to $x$ from the left and the right such that
there is a ball of radius $R$ containing both $y$ and $y'$ but not $x$, a contradiction.
\end{proof}

Just like monotone functions, continuous $R$-monotone curves do not have to be
twice differentiable, consider the arc-length parametrized version of $s(t)=(t, \vert t \vert t)$
for $t\in (-1,1)$.

Differentiable functions are not necessarily monotone. Likewise
differentiable curves are not necessarily $R$-monotone as the next
proposition shows.

\begin{proposition}
 There is a differentiable continuous finite
 curves with no cross-over (that is, $t\rightarrow s(t)$ is one-to-one) which is not
$R$-monotone for any $R>0$.
\end{proposition}
\begin{proof}Consider a curve following a inward spiral around
a fixed point such as $s(t)= (2\pi-t) (\cos t, \sin t)$ for $t\in (0,2\pi]$.\end{proof}

Functions are monotone or not, and there is no ``degree of monotonicity.''
Similarly, for curves of finite length, it simply matters whether they are
$R$-monotone for some finite $R$ since $R$-monotonicity is  scale-dependent.

\begin{proposition}
Given a $R$-monotone curve, scaling the curve by a factor
$\infty>\alpha>0$ makes it $\alpha R$-monotone.
\end{proposition}

\section{Signal Monotonicity}
\label{signalmonotonicity}
In this section, we define monotonicity for vector-valued signals or chains as a
natural extension of monotonicity for real-valued signals. We show how to compute
efficiently the degree of monotonicity.

 A scalar-valued signal (or discrete function) is monotone if and
only
if the index set of values in any closed interval $[a,b]$ is a set
of consecutive integers $[j,k]$: $p_i \in [a,b] \Leftrightarrow i\in
[j,k]$. Equivalently, the values of the signal $p_i$ never go down ($p_{i+1} \geq
p_i$) or never
go up ($p_{i+1} \leq p_i$). Another equivalent definition
is given
by the next proposition.

\begin{proposition}
A scalar-valued signal $p_i$ is monotone if and only if, for any 3
consecutive samples, $p_i, p_{i+1}, p_{i+2}$, the index set of the values
contained in any closed interval $[a,b]$ is a set of consecutive integers
$[j,k]$. Equivalently, the
index set is a convex set under an appropriate definition of convexity.
\end{proposition}

It is easy to extend this definition of monotonicity to the case of
vector-valued
signals.
Unfortunately, a straightforward generalization, based on considering the
set of indices
of the values contained in any closed ball, would lead us to conclude that
the only monotone
vector-valued signals are on straight lines and never backtrack. It is not
hard to realize
no sensible filter could turn any vector-valued signal into a monotone
signal.
In order to obtain nontrivial results, we need to restrict the class of
balls considered, as in the following definition.

\begin{definition}
A vector-valued signal $\mathbf{p}_i$ has a degree of monotonicity $R$ if
$R$ is the largest value such that, considering only 3 consecutive samples,
$\mathbf{p}_i, \mathbf{p}_{i+1}, \mathbf{p}_{i+2}$, the index set of the
values contained in any closed ball $B$ of radius at most $R$ is a
set of consecutive integers in 
$\{i,i+1,i+2\}$.
\end{definition}

\begin{figure}
\center{\includegraphics[width=0.60\columnwidth,angle=0]{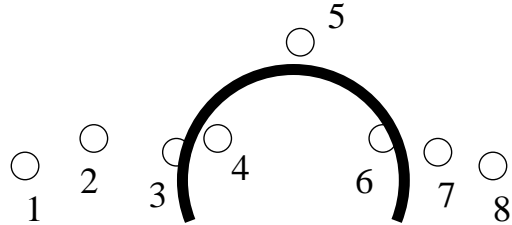}}
\caption{\label{monotonicityfailure} \small Given the chain of data points
shown, the degree of monotonicity is at most the size of the radius of the
circle given in the picture: it contains points 4 and 6 but not point 5.}
\vspace*{-0.5cm}
\end{figure}

If the signal values are on a straight line with no backtracking, then the
degree of monotonicity is $\infty$, and the degree of monotonicity is
always larger than $0$ for finite signals.
Fig.~\ref{monotonicityfailure} gives an intuitive view of the degree of
monotonicity. 
This measure of monotonicity is robust in the following
sense.

 \begin{proposition} If one point is omitted from a vector-valued signal,
the degree of monotonicity cannot decrease.
\end{proposition}

While this discrete definition is similar to the definition 
given for $R$-monotone curves, to allow efficient computation, we consider
only sets of 3 consecutive samples, thus replacing a global problem by a local problem. 
If we lift the requirement that
only 3~samples are considered, then a signal is $R$-monotone 
if and only if all subchains of
length $3$ are $R$-monotone. This suggests that
 the cost of checking global
$R$-monotonicity grows in a cubic
fashion with respect to the length of the signal which is unacceptable for most applications.

In practical applications, maximizing the degree of monotonicity $R$
leads to useful chains. For example, noise tends to reduce $R$ by creating sharp turns and 
local backtracking and a highly monotone curve ($R$ large) is more likely to be noise-free.
On the other hand, when reconstructing chains from unordered sets of points, as happens in 
computer vision, we often want to minimize sharp turns and backtracking. Therefore,
solving for the chain maximizing $R$ while passing through all available data points
is a sensible ``curve reconstruction'' strategy.

As a prerequisite to computing the degree of monotonicity, we need 
a computationally effective way to
compute the radius of the circle going through 3 points. Given
$\mathbf{p}_1, \mathbf{p}_2, \mathbf{p}_3 \in \mathbb{R}^n$, we can
compute the radius
of the circle passing
through them (denoted $\frown (\mathbf{p}_1 \mathbf{p}_2 \mathbf{p}_3)$
)
by first computing $a= \Vert \mathbf{p}_1-\mathbf{p}_2\Vert$, $b=\Vert
\mathbf{p}_2 - \mathbf{p}_3\Vert$,
$c= \Vert \mathbf{p}_1 - \mathbf{p}_3 \Vert$, $\sigma=(a+b+c)/2$, and then we
have the classical Heron's formula for the radius of the circle:
\[R^{\mathrm{outcircle}}= \frac{abc}{4 \sqrt{\sigma (\sigma-a) (\sigma-b) (\sigma-c)}}\]
whenever $a,b,c>0$. 

\begin{figure}
\center{\includegraphics[width=0.90\columnwidth,angle=0]{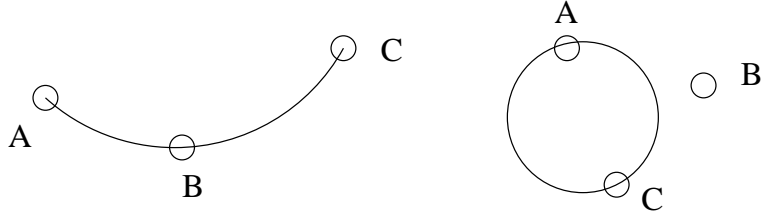}}
\caption{\label{computingR} Given a chain of 3 data points, we give two
cases: (left) the angle $\angle (ABC) > \pi/2$ so we compute the radius of
the circle going through $ABC$, otherwise (right), we compute half the
distance between $A$ and $C$.}
\vspace*{-0.5cm}
\end{figure}

The next theorem gives us a way to compute the (local) degree of monotonicity for
any 3 points, to compute the degree of monotonicity of an entire signal
simply requires, by definition, to take the \textbf{minimum} of the result for all
consecutive 3 points. The theorem essentially says that if $\angle
(\mathbf{p}_1\mathbf{p}_2\mathbf{p}_3) < \pi/2$, the degree of
monotonicity is  then half the distance between $\mathbf{p}_1$ and
$\mathbf{p}_3$, and otherwise, it is  $R^{\mathrm{outcircle}}$ (see
Fig.~\ref{computingR}). To see that this local form of monotonicity is distinct from the 
global form suggested earlier, consider a chain in the form of a figure ``\textbf{8}.'' 

\begin{theorem}\label{formulathm} The degree of monotonicity for the
sequence  $\mathbf{p}_1,\mathbf{p}_2,\mathbf{p}_3$
is
\begin{align*}
R:=\left\{\begin{array}{cl}
\frac 12c&\mbox{ if } a^2+b^2>c^2\\
R^{\mathrm{outcircle}}&\mbox{ otherwise}\end{array}
\right.\end{align*} \label{prop:formula}
where $a= \Vert \mathbf{p}_1-\mathbf{p}_2\Vert$, $b=\Vert
\mathbf{p}_2 - \mathbf{p}_3\Vert$,
$c= \Vert \mathbf{p}_1 - \mathbf{p}_3 \Vert$.
\end{theorem}
\begin{proof}
Consider the disk $B_0$ containing $\mathbf{p}_1$ and $\mathbf{p}_3$,
centered at
$(\mathbf{p}_1+\mathbf{p}_3)/2$ and having radius
$d(\mathbf{p}_1,\mathbf{p}_3)/2$. The point $\mathbf{p}_2$ is
outside the disk if and
only if $\cos\angle (\mathbf{p}_1\mathbf{p}_2\mathbf{p}_3)=\frac{a^2+b^2-c^2}{2ab}$ is positive.
Thus, $\mathbf{p}_2$ is outside the disk if and
only if  $a^2+b^2-c^2>0$. Clearly $R= \mathrm{radius}(B_0)
= d(\mathbf{p}_1,\mathbf{p}_3)/2$.
Next, suppose that $\mathbf{p}_2$ is in the disk $B_0$. We have that any
ball
containing $\mathbf{p}_1$ and $\mathbf{p}_3$ but not $\mathbf{p}_2$ must
be larger than
$\mathrm{radius}(B_0)$  since $B_0$ is the smallest ball containing both
$\mathbf{p}_1$ and $\mathbf{p}_3$. Now, suppose there is a (closed) ball
of minimal radius
$R$ containing $\mathbf{p}_1$ and $\mathbf{p}_3$, but not $\mathbf{p}_2$.
This implies a non-zero
distance, $\delta >0$, between $B$ and $\mathbf{p}_2$. We have that the
center of
the ball has to be away from the line formed by
$\mathbf{p}_1,\mathbf{p}_3$: if not then it must
be a ball containing $B_0$. This means we can move the center of the ball
slightly closer to $\mathbf{p}_1$ and $\mathbf{p}_3$ while reducing the
radius just enough
so that $\mathbf{p}_2$ remains outside the ball. By repeating this
process, we show
that $\delta=0$, a contradiction. Hence, there is no (closed) ball of
minimal radius $R$ containing $\mathbf{p}_1$ and $\mathbf{p}_3$, but not
$\mathbf{p}_2$. Hence
$\mathbf{p}_1,\mathbf{p}_2,\mathbf{p}_3$ have a degree of monotonicity
$R^{\mbox{{\rm outcircle}}}$.
\end{proof}

 \section{Monotonicity and Sphere-Preserving Filters }
 \label{sec:spherepreservingsection}

 In this section, we propose a SP filter which  never decreases the degree
of monotonicity of the signal. Given a signal $\mathbf{p}_i$, we consider
recursive (IIR) filters of the form
 \[\mathbf{p}'_i= f(\mathbf{p}'_{i-2},\mathbf{p}'_{i-1},\mathbf{p}_{i},
\mathbf{p}_{i+1}, \mathbf{p}_{i+1}).\]
 To ease the notation, we write $A= \mathbf{p}'_{i-2}$,
$B=\mathbf{p}'_{i-1}$, $X=\mathbf{p}_{i}$, $C=\mathbf{p}_{i+1}$,
$D=\mathbf{p}_{i+1})$ so that the equation becomes $X'=f(A,B,X,C,D).$
 Let $R(A,B,X,C,D)$ be the degree of monotonicity of $A,B,X,C,D$ computed
as $\min(R(A,B,X),R(B,X,C),R(X,C,D)).$
 The following proposition gives us a condition of $f$ to increase the
monotonicity of a vector-valued signal.
 \begin{proposition}
 Given $X'=f(A,B,X,C,D)$, if $f$ is such that the degree of monotonicity
 \[R(A,B,X',C,D) \geq R(A,B,X,C,D),\] then the recursive filter
 \[\mathbf{p}'_i= f(\mathbf{p}'_{i-2},\mathbf{p}'_{i-1},\mathbf{p}_{i},
\mathbf{p}_{i+1}, \mathbf{p}_{i+1})\]
 never decreases the degree of monotonicity of a signal.
 \end{proposition}
 It seems that $f$ should be chosen so that $R(A,B,X',C,D)$ is as large as
possible. To maximizes $R(A,B,X',C,D)$
with $X'=f(A,B,X,C,D)$, $f$ should be either $B$
or $C$. In other words, we improve monotonicity best when we make the
sample $X$ ``virtually disappear.''
 \begin{proposition} $R(A,B,X',C,D)$ is minimized when $X'=B$ or $X'=C$
and these choices are unique unless $\frown (ABC) = \frown (BCD)$ in
which case any point on the arc of the circle between $B$ and $C$
inclusively qualifies.
 \end{proposition}
 Fortunately, we can easily define a more interesting SP filter. Given an
arc of a circle, denoted $\alpha$, and a point $X$, we can project $X$ on
$\alpha$ by solving for the point closest $X$ in $\alpha$.
 The projection onto a circle 
can be determined easily
using only linear algebra~\cite{Kalman2002}. In the plane,
start with equation
$(x-r_1)^2+(y-r_2)=\rho^2$ and substitute 3 values of $x,y$,
getting 3~equations. By pairwise subtraction, we can remove the unknown $\rho^2$, and 
be left with linear system having 2~equations and 2~unknowns (the center of the circle).
We apply this by first projecting on the circle
and if the projected point does not belong to the given arc we move it to
the closest point on the arc (an endpoint of the arc).
 Let us define
 $X_1$ to be the projection of $X$ on the arc $BC$ of the circle $ABC$, and
define $X_2$ to be the projection of $X$ on the arc $BC$ of the circle
$BCD$. Intuitively, either point $X_1$ or $X_2$ would make a good choice
for $X'$ 
.
 To ensure that the degree of monotonicity is never decreased, we set
 \[f(A,B,X,C,D)= \arg \max_{X'\in \{X,X_1,X_2\}} R(A,B,X,C,D).\]
  This function can be computed quickly and is sphere-preserving.
%
%

 \section{Experimental Results}

\begin{figure}
\centering{\includegraphics[height=0.7\columnwidth,angle=270]{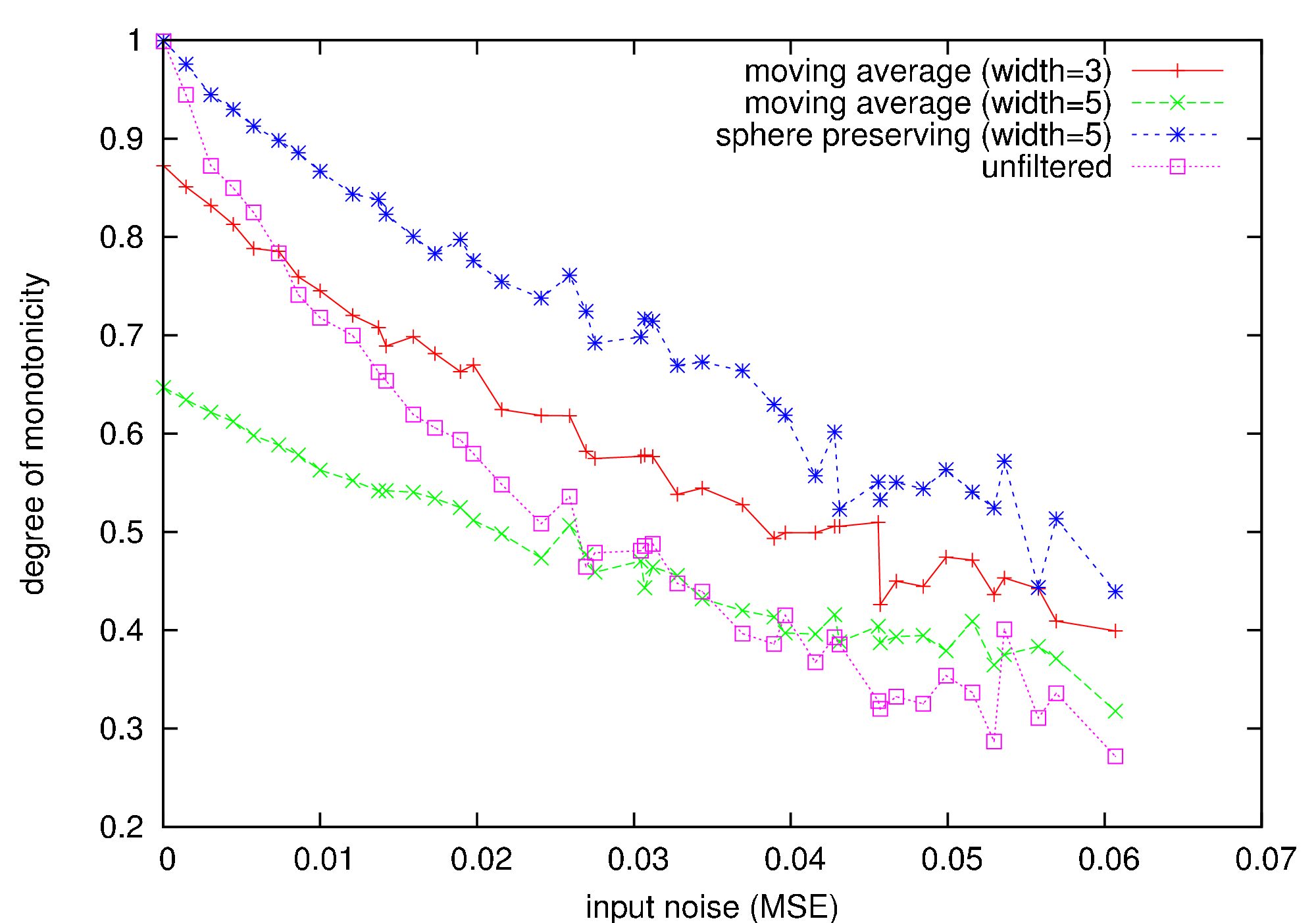}
}
\caption{\label{monversusnoiseversuswidth}\small The degree of
monotonicity versus the absolute input noise level (MSE) over a synthetic
data set generated from points on a unit circle. The SP filter outperforms
MA when noise levels are low. }
\vspace*{-0.5cm}
\end{figure}
\begin{figure}
\centering{\includegraphics[width=0.7\columnwidth,angle=0]{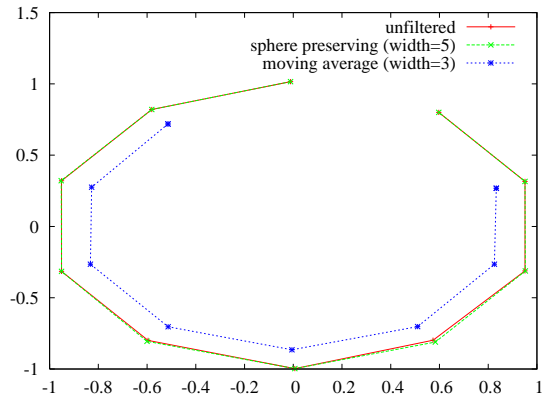}
}
\caption{\label{visualcompared}\small Visual comparison of the SP
filter with the MA filter for low noise levels and coarse sampling. }
\vspace*{-0.5cm}
\end{figure}

 We generate a chain in the $xy$ plane by regularly sampling a unit circle
3 times for a total of 30~samples. A MA filter with window width $k$ averages each $k$ consecutive
data points. We add white noise to every point in
the chain and we filter it using simple MA filters with window widths of 3
and 5 samples as well as with the SP filter of the previous section. Each
test is repeated 10~times and we keep only the averages.
 Fig.~\ref{monversusnoiseversuswidth} shows the degree of monotonicity
versus the noise level (Mean Square Error) with the three smoothing
filters and the unfiltered chain.
The noise level ranges from none to over 0.05 (MSE) which corresponds
roughly to a 5\% noise-to-signal ratio. An example of filtering is given in Fig.~\ref{visualcompared}.

 \section{Discussion}


In the unfiltered chain, the degree of monotonicity is inversely
correlated with the noise level: the Pearson correlation is $p=-0.95$
(90\%).  The degree of monotonicity seems a good indicator of noise, which
in particular suggests that a method for increasing the degree of
monotonicity would also function as a good noise reduction technique.  As
required, the SP filter always increases the degree of monotonicity with
respect to the unfiltered data. Simple MA filters \textbf{decrease the
degree of monotonicity} when noise levels are low, and more aggressive
filtering (window width of 5 versus 3) even more so. The result of
aggressive lowpass filtering on the curvature of a chain is explained by
Fig.~\ref{filteredcircle}. The relative performance of filters over chains
can vary depending on the level of noise and  the distance between the
points: as noise levels increase, the SP filter is less competitive.
The design of sphere-preserving filters optimally increasing the degree of
monotonicity is an open problem.



\vfill{}

\begin{address}
\begin{minipage}[b]{0.47\linewidth}
Dan Kucerovsky
University of New Brunswick
Fredericton NB CANADA
{\tt dan@math.unb.ca
\url{www.math.unb.ca/~dan/}}
\end{minipage} \hfill\begin{minipage}[b]{0.47\linewidth}
Daniel Lemire
University of Quebec at Montreal
Montreal QC CANADA
{\tt lemire@acm.org
\url{www.daniel-lemire.com/}}
\end{minipage}

\end{address}

\end{document}